\documentclass{article}
\usepackage{spconf,amsmath,amsfonts,amssymb,amsthm, graphicx}
\usepackage{multirow}
\usepackage{rotating}
\usepackage{url}
\usepackage{color}
\usepackage{enumerate}

\def\R{\mathbb{R}}

\makeatletter
\def\captionof#1#2{{\def\@captype{#1}#2}}
\makeatother

\title{Variable density sampling based on physically plausible gradient waveform. Application to 3D MRI angiography.}
\name{\!\! Nicolas Chauffert$^{\small (1,2)}$,\! Pierre Weiss$^{\small(3)}$,\! Marianne Boucher$^{\small(1)}$,\! S\'ebastien M\'eriaux$^{\small(1)}$\! and\! Philippe Ciuciu$^{\small(1,2)}$}
\address{$^{\small (1)}$CEA/DSV/I$^2$BM NeuroSpin center, B\^at. 145, F-91191 Gif-sur-Yvette, France\\
$^{\small (2)}$ INRIA Saclay Ile-de-France, Parietal team, 91893 Orsay, France.\\
$^{\small (3)}$ PRIMO Team, ITAV, USR 3505, Universit\'e de Toulouse.}

\begin{document}
\ninept
\maketitle
\begin{abstract}
Performing $k$-space variable density sampling is a popular way of reducing scanning time in Magnetic Resonance Imaging~(MRI). Unfortunately, given a sampling trajectory, it is not clear how to traverse it using gradient waveforms. In this paper, we actually show that existing methods~\cite{Lustig08,Vaziri13} can yield large traversal time if the trajectory contains high curvature areas. Therefore, we consider here a new method for gradient waveform design which is based on the projection of unrealistic initial trajectory onto the set of hardware constraints. Next, we show on realistic simulations that this algorithm allows implementing variable density trajectories resulting from the piecewise linear solution of the Travelling Salesman Problem in a reasonable time. Finally, we demonstrate the application of this approach to 2D MRI reconstruction and 3D angiography in the mouse brain.
\end{abstract}
\begin{keywords}
MRI, Compressive sensing, Variable density sampling, gradient waveform design, hardware constraints, angiography.
\end{keywords}
\section{Introduction}
\label{sec:intro}
\begin{figure}[!h]
\begin{center}
\begin{tabular}{cc}
(a)&(b)\\
\includegraphics[width=.35\linewidth]{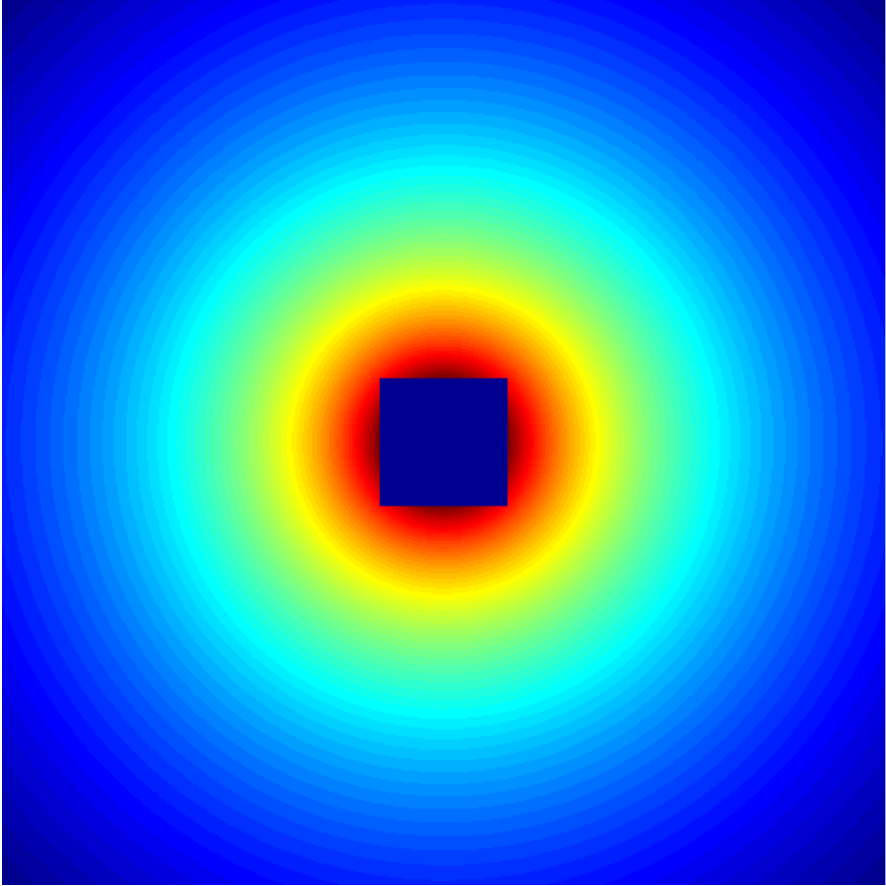}&
\includegraphics[width=.35\linewidth]{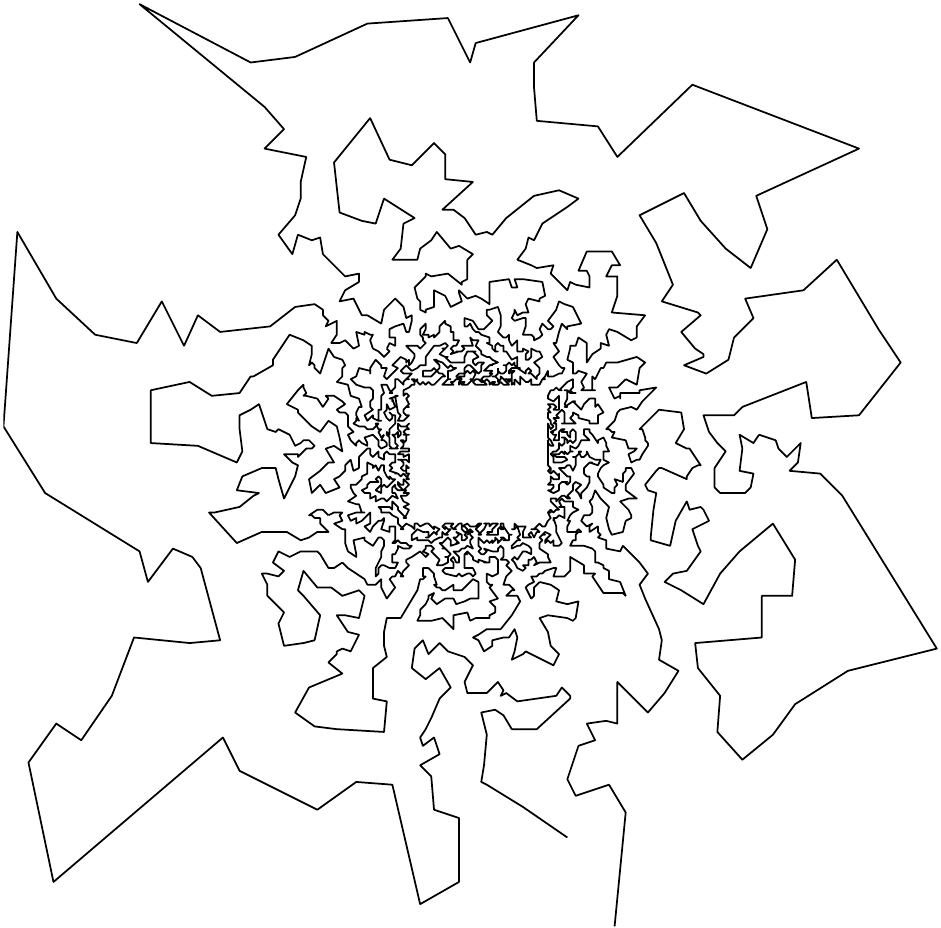}
\end{tabular}
\end{center}
\caption{\label{fig:introVDS} Example of 2D Variable Density Sampler. (a):~$\pi(k) \propto 1/|k|^2$ as advocated in~\cite{Chauffert14}. (b):~TSP-based sampling trajectory.\vspace{-.05 \linewidth}}
\end{figure}

Compressed Sensing~(CS) provides a theoretical framework to justify the downsampling of $k$-space~(2D or 3D Fourier domain) in Magnetic Resonance Imaging~(MRI). CS-MRI is usually based on independent random drawing of $k$-space locations according to a prescribed density. From recent theoretical works~\cite{Candes11,Rauhut10}, one can derive an optimal sampling density $\pi$ that reduces at most the number of samples collected in MRI without degrading the image quality at the reconstruction step~\cite{Chauffert13,Puy11}. In~\cite{Chauffert14}, simulations show that distributions with radial decay (see Fig.~\ref{fig:introVDS}(a)) with full $k$-space center acquisition perform better in numerical experiments. 

However, such sampling schemes are not performed along continuous lines and thus not physically plausible in MRI because of the constraints involved on the magnetic field gradient~(magnitude and slew-rate). In~\cite{Chauffert13b}, we have proposed a new approach to design continuous sampling trajectories based on the solution of Travelling Salesman Problem~(TSP), as illustrated in Fig.~\ref{fig:introVDS}(b). The specificity of this approach is that the empirical distribution of the trajectory can approximate any prescribed distribution $\pi$. Such a curve is called a $\pi$-Variable Density Sampler~($\pi$-VDS). Unfortunately, continuity of the sampling trajectory is not a sufficient condition in MRI and it is not clear how to design admissible gradient waveforms to traverse such a trajectory.

To the best of our knowledge, the most efficient gradient design strategies consist of finding an admissible parameterization of a given curve by using optimal control~\cite{Lustig08,Vaziri13}. However, the main drawback of this approach is that the traversal of high curvature parts of the curve is rather slow. In the extreme case of angular points such as in TSP-based trajectories, it leads to extremely large acquisition time. In~\cite{Chauffert14b}, we have proposed a new method to design magnetic field gradients by projecting any parameterized curve onto the set of hardware constraints. This method allows one to change the curve support and thus to yield faster traversal time.

The goal of this paper is then to prove that this projection algorithm enables to implement a TSP-based VDS on MRI scanners in reasonable scanning times, in contrast to optimal reparamerization. The first part is dedicated to summarizing the projection strategy introduced in~\cite{Chauffert14b}. Next, we provide in Section~\ref{sec:simul} a proof-of-concept on retrospective CS simulations. In particular, we show that our algorithm yields faster sampling trajectories than the state-of-the-art for a given image reconstruction quality, and alternatively that if the scanning time is fixed, our method delivers improved reconstructions. Finally, our strategy~(TSP-based sampling + projection algorithm) is applied to 3D angiography of the mouse brain before and after contrast agent injection, to demonstrate its efficiency in terms of scanning time reduction while preserving the recovery of small structures such as blood vessels.

\section{Gradient waveform design using a projection algorithm}
In this section, we recall the hardware constraints as generally model\-led in MRI and describe methods for designing gradient waveforms in order to traverse $k$-space sampling curves.

The gradient waveform associated with a curve $s$ is defined by $ g(t)=\gamma^{-1} \dot{s} (t)$, where $\gamma$ denotes the gyro-magnetic ratio~\cite{Hargreaves04}.
The gradient waveform is obtained by energizing gradient coils (arrangements of wire) with electric currents.
Owing to obvious physical constraints, these electric currents have a bounded amplitude and cannot vary too rapidly~(slew rate). Mathematically, these constraints read: 
\begin{align*}
\|g(t)\| \leqslant G_{\max} \qquad \mbox{and} \qquad \|\dot g(t)\| \leqslant S_{\max}, \quad \forall t \in [0,T].
\end{align*}
A sampling trajectory $s:[0,T]\to \R^d$ will be said \emph{admissible} if it belongs to the set:
\begin{align*}
\mathcal{S}:=\left\{s\in \left(\mathcal{C}^2([0,T])\right)^d, \|\dot s(t)\| \leqslant \alpha , \|\ddot s(t)\| \leqslant \beta, \quad \forall t \in [0,T] \right\}.
\end{align*}
In this paper, we limit ourselves to the so-called  \emph{rotation-invariant} constraints. Some hardware systems, where the coils are energized independently, enable considering \emph{rotation-variant} constraints. In this case, the $\ell_2$ norm $\| \cdot \|$ is replaced by $\ell_\infty$ norm. The differences are discussed in~\cite{Vaziri13}, but the two compared methods~(optimal control and projection algorithm) are able to handle both kinds of constraints. 

\subsection{State-of-the-art}
The question of finding jointly an accurate trajectory and the admissible gradient waveform to traverse it is a difficult issue that has received special attention in~\cite{Dale04,Mir04}. The most classical approaches consist of fixing a curve support~\cite{Lustig08,Vaziri13}, or control points~\cite{Hargreaves04,Davids14b}, and finding an admissible parameterization afterwards.

In particular, the most popular approach to design an admissible curve assumes the knowledge of a parameterized curve $c:[0,T]\to \R^d$ and consists of finding its optimal reparameterization by using optimal control~\cite{Lustig08,Vaziri13}. In other words, it amounts to finding a reparameterization $p$ such that $s=c\circ p$ satisfies the physical constraints while minimizing the acquisition time. This problem can be cast as follows:
\begin{align}\label{eq:reparam}
T_{\texttt{OC}}=\min T' \quad \text{such that} \quad \exists\, p:[0,T'] \mapsto [0,T],\; c \circ p \in \mathcal{S}.
\end{align}
The resulting solution $s=c\circ p$ has the same support as $c$, which might be an important feature in some applications.
The problems of this approach are: i) there is no control of the sampling density, especially in the high curvature parts of $c$ where samples tend to agglutinate and ii) there is no control over the total sampling time $T_{\texttt{OC}}$, which can be large if the trajectory contains singular points for instance (e.g., see Fig.~\ref{fig:introVDS}(b)). These two drawbacks are illustrated in~\cite{Chauffert14b}.

The next part is dedicated to introducing an alternative method relaxing the constraint of the curve support.

\subsection{Projection onto the set of constraints}

The idea introduced in~\cite{Chauffert14b} is to find the projection of the input curve $c:[0,T]\to\R^d$ onto the set of admissible curves $\mathcal{S}$:

\begin{align}
\label{pb:primal}
s^*:= \underset{s\in \mathcal{S}}{\operatorname{argmin}} \frac{1}{2} \|s-c\|_2^2.
\end{align}
where $\|s-c\|_2^2:=\int_{t=0}^T \|s(t)-c(t)\|_2^2 dt$. For the sake of conciseness, the theoretical grounds for the resolution and the key properties discussed below are given in~\cite{Chauffert14b}.

\noindent\textbf{Resolution.} Problem~\ref{pb:primal} consists of minimizing a convex smooth function over a convex set. In~\cite{Chauffert14b}, we have proposed a fast iterative algorithm exploiting the structure of the dual formulation of~\eqref{pb:primal}, which can be solved using proximal gradient methods~\cite{Combettes11b}.\\
\noindent\textbf{Key properties.}
The two main advantages of this method are that i) the sampling density of $s^*$ is close to the sampling density of $c$. In particular, if $c$ is a $\pi$-VDS, the sampling density of $s^*$ is close to $\pi$\footnote{the \emph{closeness} is quantified by Wasserstein transportation distance $W_2$, see~\cite{Chauffert14b} for details.}; and ii) the acquisition time $T$ is fixed and equal to that of the input curve $c$. In particular, the time to traverse a curve is in general shorter than with optimal reparameterization ($T<T_{\texttt{OC}}$).

In the next part, we will emphasize that our algorithm enables to traverse VDS curves as depicted in Fig.~\ref{fig:introVDS}(b) in a reasonable time, unlike optimal control-based reparemeterizations.

\section{CS-MRI simulations}
\label{sec:simul}

In this part, we compare the time to traverse $k$-space along different trajectories using gradients computed either by the standard optimal control approach or by our proposed projection algorithm. For comparison between sampling schemes, we work on \emph{retrospective} CS, meaning that a full dataset has been acquired, and then a posteriori downsampling is performed. \noindent We compare the reconstruction results in terms of peak signal-to-noise ratio~(PSNR) with respect to the acquisition time and to the ``acceleration factor''\footnote{$r$ quantifies the reduction of the number of measurements $m$. If the $k$-space is a grid of $N$ pixels $r:=N/m$ is commonly used in CS-MRI.} $r$.

\subsection{Experimental framework}

\noindent\textbf{Data acquisition.}
The initial experimental setup aimed at observing blood vessels of living mice using an intraveinous injection of an iron oxide-based contrast agent~(Magnetovibrio Blakemorei MV1). Because of natural elimination, it is necessary to speed up acquisition to improve contrast and make easier post-processing such as angiography.
The experiments have been performed on a 17.2T preclinical scanner which physical rotation-invariant constraints are, for all $t \in [0,T]$: 
\begin{align*}
\|g(t)\| \leqslant 1 \mbox{ T.m$^{-1}$} \qquad \mbox{and} \qquad \|\dot g(t)\| \leqslant 5.3 \mbox{ T.m$^{-1}$.ms$^{-1}$} .
\end{align*}
A FLASH sequence (Fast Low Angle SHot) has been used to reveal the $\text{T}^*_2$ contrast induced by the injection of the contrast agent (TE/TR = 8/680~ms). The sequence was repeated 12 times to improve the signal-to-noise ratio~(SNR), leading to a total acquistion time of 30 minutes to acquire the $k$-space slice by slice. The spatial resolution achieved is $90\!\times\! 90\! \times \! 180 \ \mu\text{m}^3$.

\noindent\textbf{Hypothesis.}
\label{part:hypothesis}
The aim of this paper is to prove that one can expect a large acquisition time reduction using partial $k$-space measurements. The time to traverse a sampling curve is computed satisfying the gradient constraints. To achieve a fair comparison, let us mention the additional hypothesis that our acquisitions are single-shot, meaning that the partial $k$-space is acquired after a single RF pulse. We did not take echo and repetition times into account to ensure the recovery of a $\text{T}^*_2$-weigthed image. We only compare the time to traverse a curve using the gradients with their maximal intensity. We assume that there is no error on the $k$-space sample locations. In practice we have to measure the three magnetic field gradients that are actually played out by the scanner to correct the trajectory and avoid distortions.
We shall work on a discrete cartesian $k$-space, and consider that a sample is measured if the sampling trajectory crosses the corresponding cell of the $k$-space grid. 
Using this hypothesis, the estimated time to visit the 2D $k$-space is 110~ms.

\noindent\textbf{Strategy.}
We used the TPS-based sampling method~\cite{Chauffert14} as input of our projection algorithm~(see Fig.~\ref{fig:2dcompar}(b)), since it is a way of designing sampling trajectories that match any sampling density $\pi$. The latter is central in CS-MRI since it impacts the number of required measurements~\cite{Adcock13,Krahmer12,Chauffert14}. To compare our projection method to existing reparameterization, the proposed sampling strategy is: \\
\noindent\textbf{(i)} Sample deterministically the $k$-space center as adviced in~\cite{Adcock13,Chauffert13,Chauffert14}, using an EPI sequence~(see Fig.~\ref{fig:2dcompar}(a)). The scanning time can be estimated to 12~ms in 2D using optimal control.\\
\textbf{(ii)} Select a density $\pi$ proportional to $1/|k|^2$ as mentioned in~\cite{Krahmer12,Chauffert14}. Draw independently points according to $\pi^{\frac{d-1}{d}}$ and join them by the shortest path to form a $\pi$-VDS~\cite{Chauffert14}.\\
\textbf{(iii)} Parameterize the TSP path at \emph{constant speed} and project this parameteri\-zation onto the set of gradient constraints, or
\textbf{(iii~bis)} Parameterize the TSP path using optimal control~(the exact solution can be computed explicitely).\\
\textbf{(iv)} Form the sampling curve, define a set $\Omega$ of the selected samples, mask the $k$-space with $\Omega$, and reconstruct an image using $\ell_1$ minimization of the constrained problem. Let $\mathbf{F}^*$ denote the $d$-dimensional discrete Fourier transform and $\mathbf{F}^*_\Omega$ the matrix composed of the lines corresponding to $\Omega$. Denote also by $\mathbf{\Phi}$ an inverse $d$-dimensional wavelet transform (here a Symmlet transform). Then the reconstructed image is the solution of the problem: 
\begin{align}
\label{eq:MinProblem}
x^* = \underset{y=\mathbf{F}^*_\Omega x}{\operatorname{Argmin}} \left\|\mathbf{\Phi}^{-1} x\right\|_1
\end{align}
An approximation of $x^*$ is computed using Douglas-Rachford algorithm~\cite{Combettes11b}. Solving the penalized form associated with~\eqref{eq:MinProblem} might be addressed by competing algorithms~(ADMM, 3MG); see~\cite{Florescu14} for a recent comparison. The reconstruction results could be improved by resorting to non-Cartesian reconstruction~\cite{Keiner09}, which would avoid the approximation related to the projection onto the $k$-space grid.

\subsection{Results}

\subsubsection{2D reconstructions}
\begin{figure*}[!ht]
\begin{center}
\begin{tabular}{ccccc}
 (a) &(b)&(c)&(d)&(e) \\
\rotatebox{90}{\hspace{.01 \linewidth}  Sampling schemes} 
\includegraphics[width=.17\linewidth]{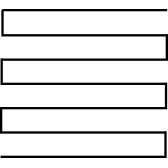} &
\includegraphics[width=.17\linewidth]{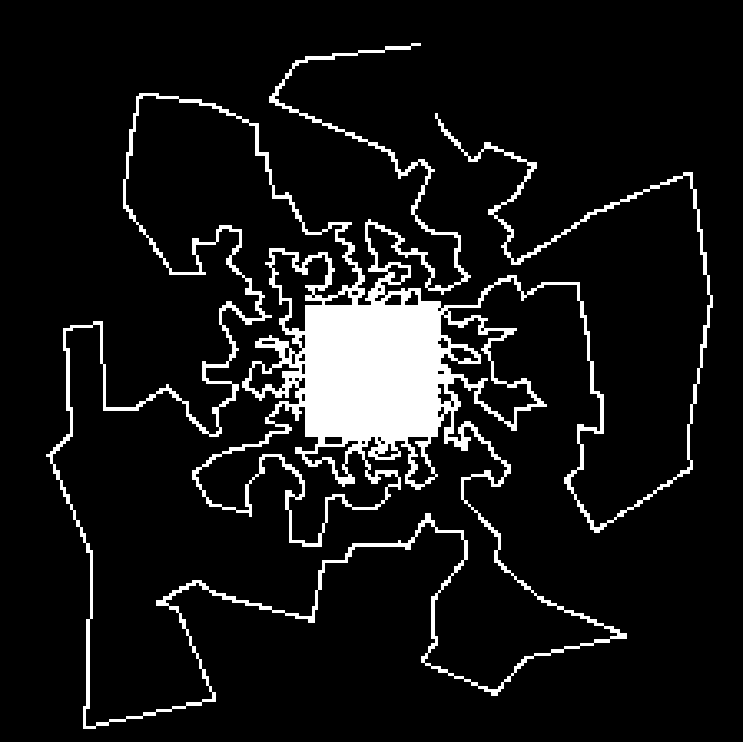}&
\includegraphics[width=.17\linewidth]{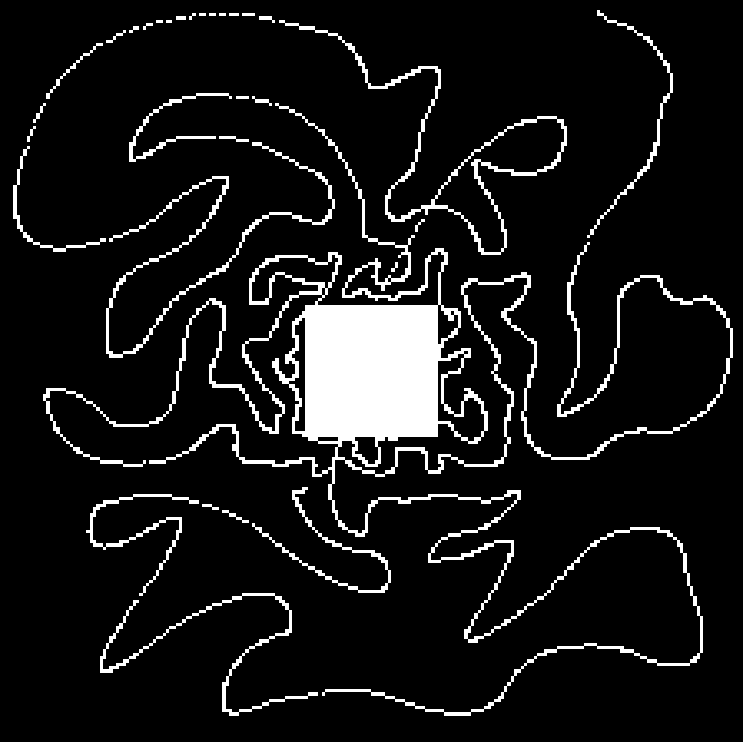}&
\includegraphics[width=.17\linewidth]{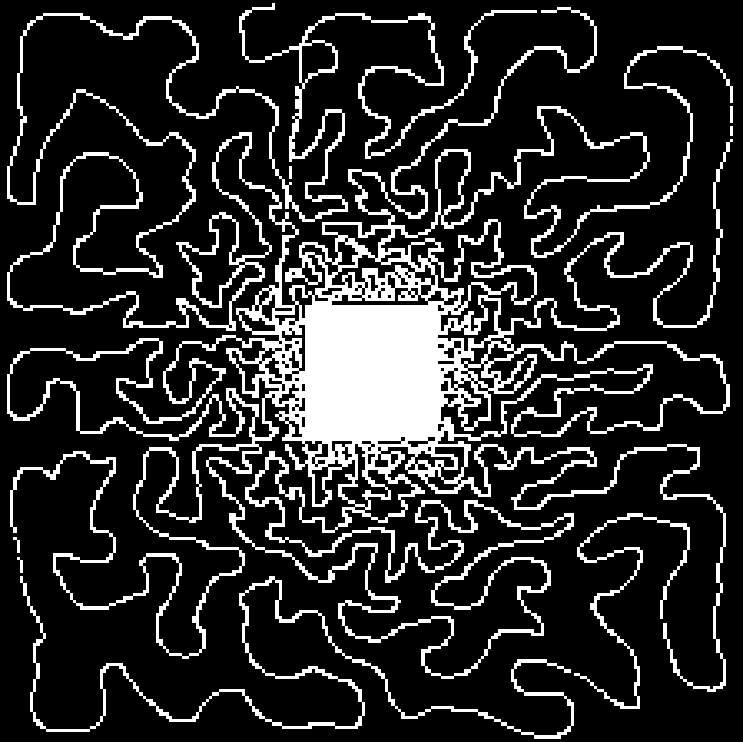}&
\includegraphics[width=.17\linewidth]{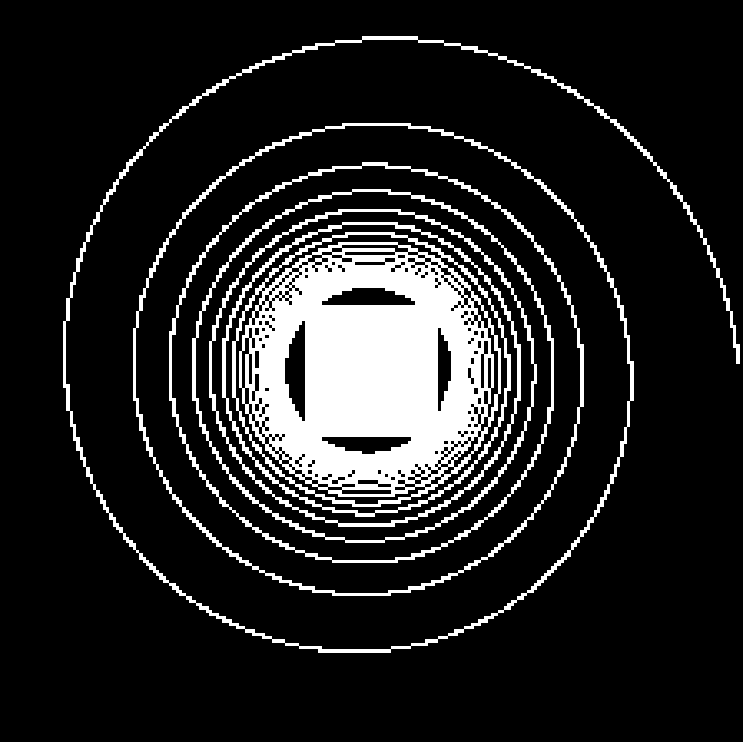}
\\
 {\scriptsize $T_{\texttt{OC}}=110$ ms ($r=1$)}& {\scriptsize $T_{\texttt{OC}}=62$ ms ($r=11.2$)} & {\scriptsize $T=30$ ms ($r=11.2$)} &
 {\scriptsize $T=62$ ms ($r=6.6$)} & {\scriptsize $T=118$ ms ($r=7.5$)}\\
   (f) & (g) &(h) &(i) & (j)\\
\rotatebox{90}{ \hspace{.015 \linewidth}Reconstructed slices}
\includegraphics[width=.17\linewidth]{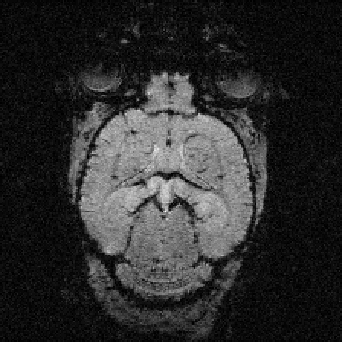}&
\includegraphics[width=.17\linewidth]{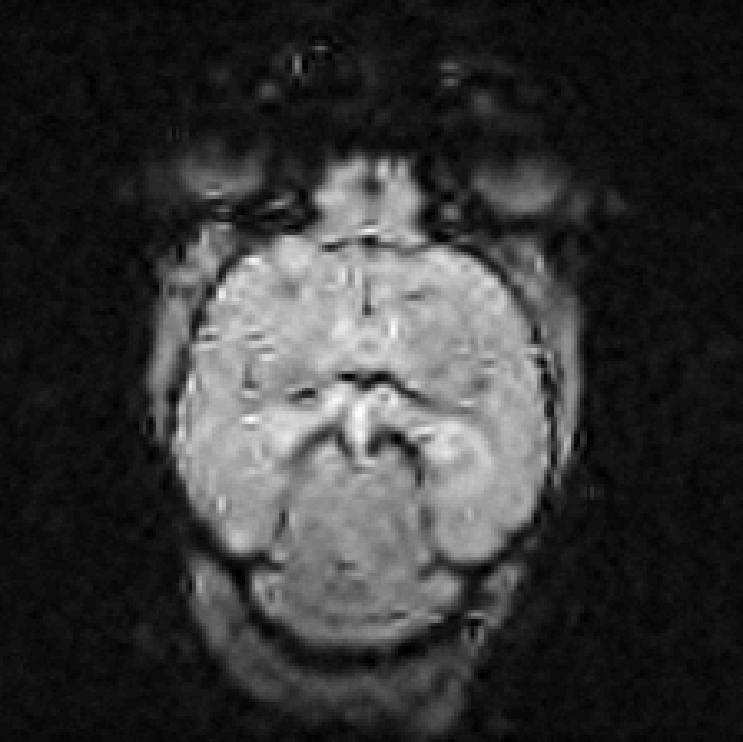}&
\includegraphics[width=.17\linewidth]{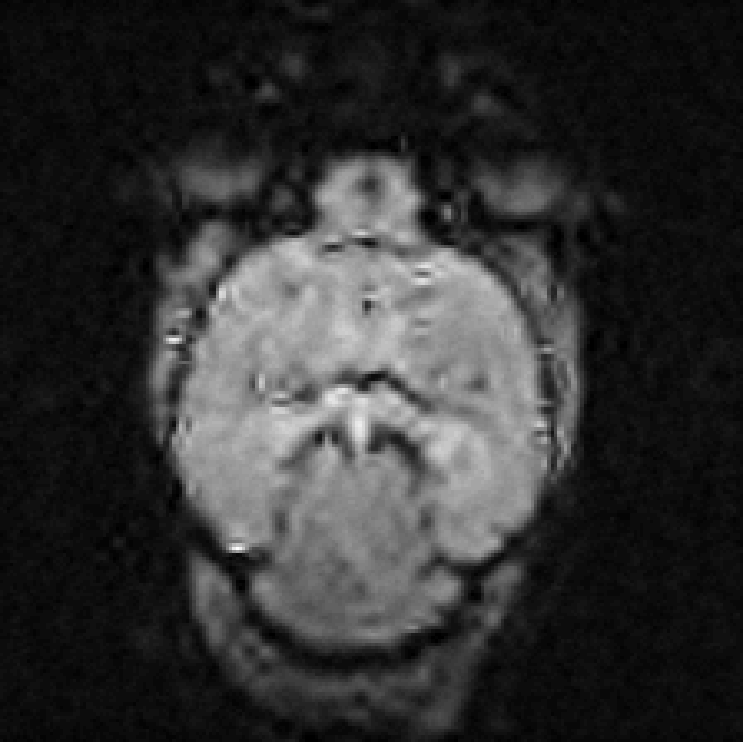}&
\includegraphics[width=.17\linewidth]{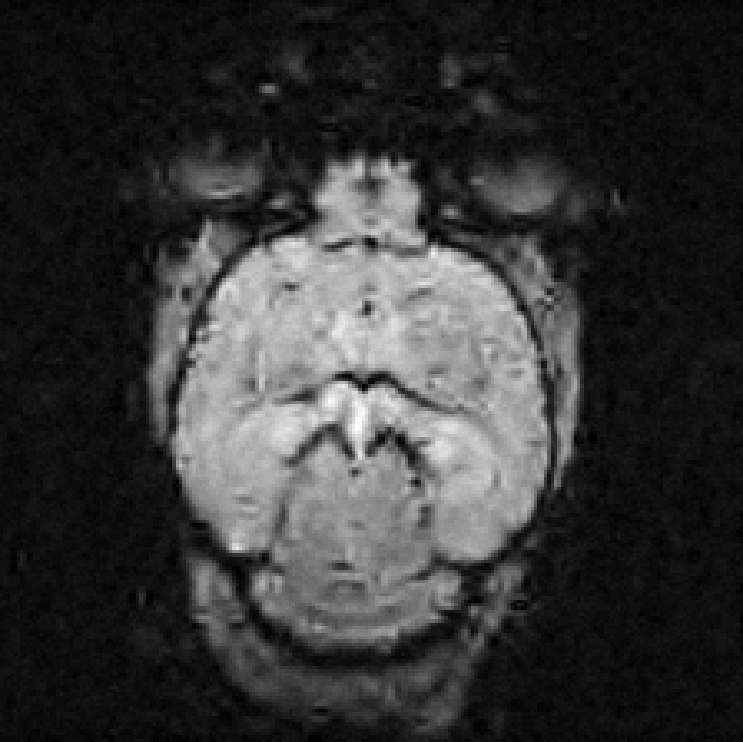}&
\includegraphics[width=.17\linewidth]{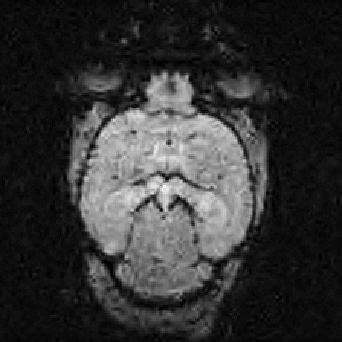}\\
\scriptsize{Reference} & \scriptsize{PSNR = 25.9 dB} & \scriptsize{PSNR = 25.5 dB}  & \scriptsize{PSNR = 26.9 dB} & \scriptsize{PSNR = 26.8 dB} \\[-.03\linewidth]
\end{tabular}
\end{center}
\caption{\label{fig:2dcompar} Full $k$-space acquisition with an EPI sequence~(a) and corresponding reference image~(f). Comparison between an exact parameterization of the TSP trajectory (b) and projection from TSP trajectory onto the set of constraints (c),(d). In experiments (b,c), the number of measured locations is fixed to 9\% ($r=11.2)$, whereas in (b,d), the time to traverse the curve is fixed to 62~ms.
(e): Spiral trajectory with full acquisition of the $k$-space center. (g-j): Reconstructed images corresponding to sampling strategies~(b-e) by solving Eq.~\eqref{eq:MinProblem}.\vspace*{-.025\linewidth}} 
\end{figure*}

\vspace*{-.02 \linewidth}
In this experiment, we considered a 2D $k$-space~($d=2$) corresponding to an axial slice. We considered five sampling strategies, depicted in Fig.~\ref{fig:2dcompar}(first row): a classical EPI coverage used as reference~(a); a TSP-based sampling trajectory parameterized using optimal control~(b); two projected TSP-based trajectories, one with the same number of samples collected as in (b)~($r=11.2$)~(c) and the other with the same scanning time as in (b)~(62~ms)~(d); a variable density spiral trajectory for comparison purpose in terms of time and sampling ratio~(e).

As expected, the reconstruction results shown in Fig.~\ref{fig:2dcompar}(g,h) are really close, since the number of collected samples is the same, and the sampling densities are similar. However, in this comparison the gain in traversal time is significant (one half). In contrast, the longer and smoothed TSP depicted in Fig.~\ref{fig:2dcompar}(d) allows us to improve image reconstruction (1~dB gain) as illustrated by~Fig.~\ref{fig:2dcompar}(i) while keeping the same acquisition time as in Fig.~\ref{fig:2dcompar}(b). 
For comparison purposes, we implemented spiral acquisition which consists of replacing steps~(ii)-(iii) in the above mentioned sampling strategy by a spiral with density proportional to $1/|k|^2$, projected onto the set of constraints. This strategy doubles the acquisition time~(118 ms compared to 62 ms) whereas the acceleration factor was larger~($r=7.5$ vs. $r=6.6$). In this experimental context (regridding and variable density spiral), the spiral is not appealing compared to EPI acquisition, since it is time consuming and degrades the image quality.

In each of these reconstructions, the major vessels can be recovered, although the smallest ones can only be seen for $r<8$. Finally, the best compromise between acquisition time and reconstruction quality is achieved using the specific combination of TSP-based sampling and our projection algorithm onto the set of constraints shown in Fig.~\ref{fig:2dcompar}(d).

\subsubsection{3D angiography}
\begin{figure}[!ht]
\hspace{-.05\linewidth}
\begin{tabular}{ccc}
 Original ($T=3.53$s) &
\hspace{-.07\linewidth} $r=7.3$ ($T=3.15$s) & 
\hspace{-.07\linewidth} $r=17.4$ ($T=0.88$s)\\
(a)&(b)&(c) \\
\rotatebox{90}{\hspace{.08\linewidth} pre-injection}
\includegraphics[width=.28\linewidth]{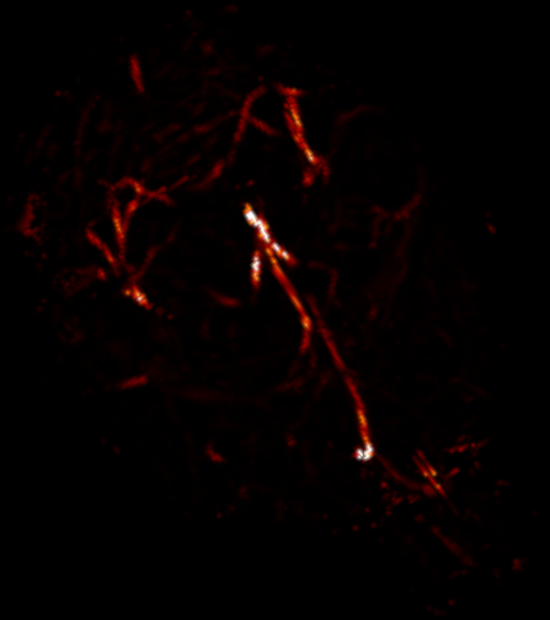}&
\includegraphics[width=.28\linewidth]{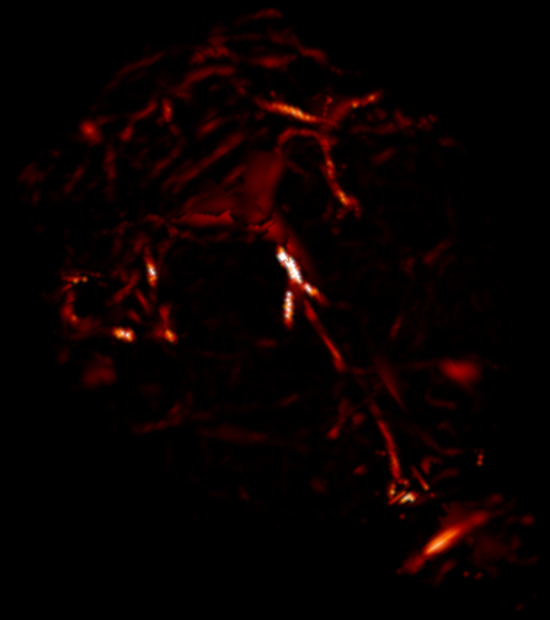}&
\includegraphics[width=.28\linewidth]{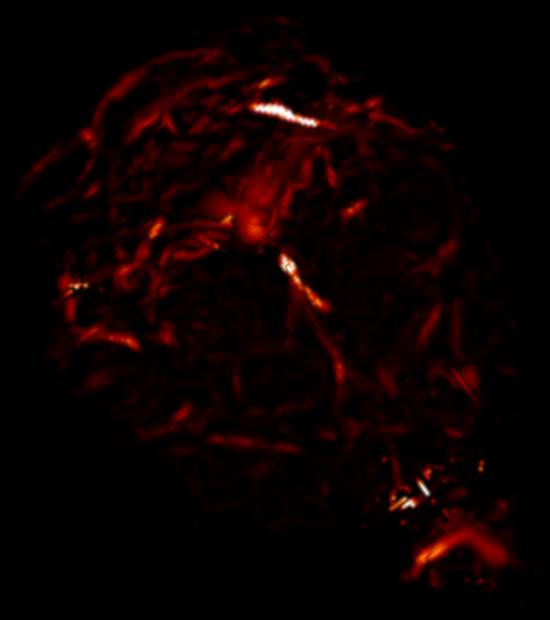}\\
  & PSNR=29.0 dB & PSNR=26.6 dB\\
 (d)&(e)&(f) \\
\rotatebox{90}{\hspace{.25\linewidth} post-injection} 
\includegraphics[width=.28\linewidth]{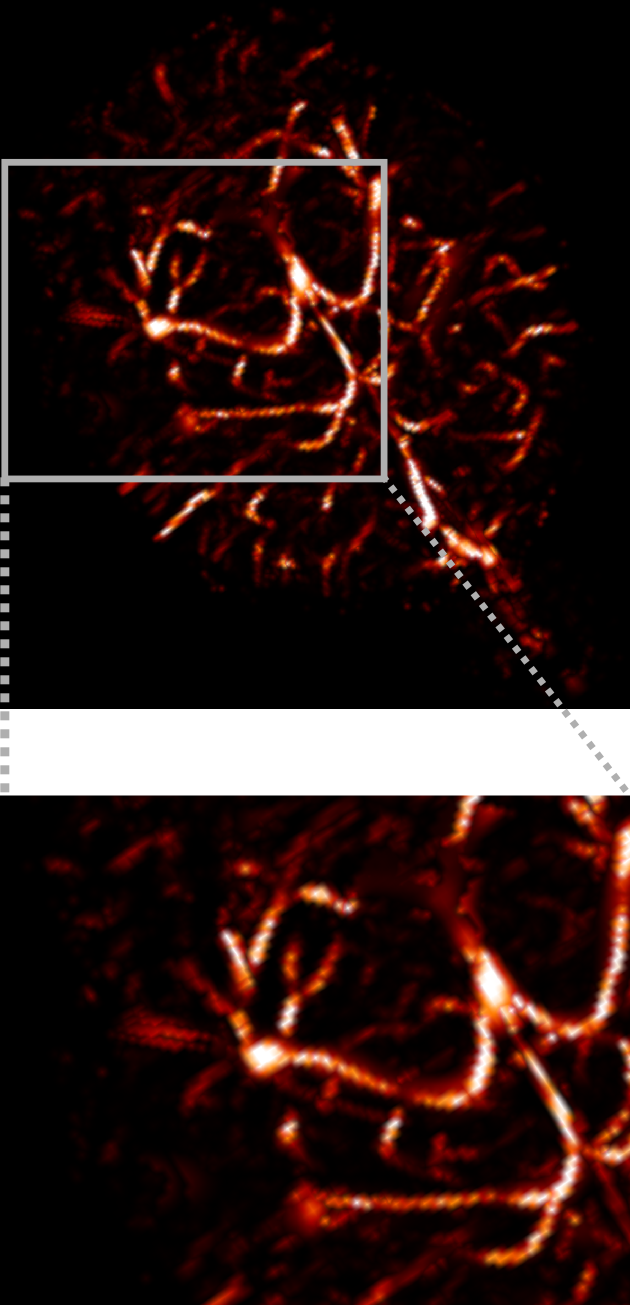}&
\includegraphics[width=.28\linewidth]{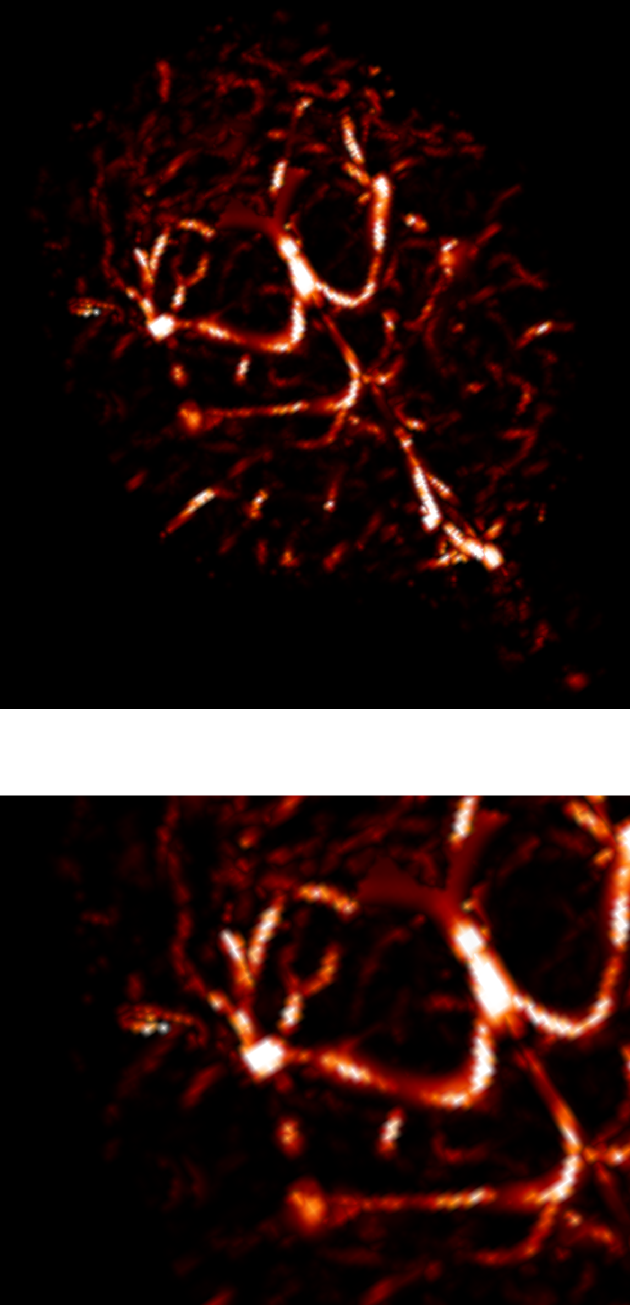}&
\includegraphics[width=.28\linewidth]{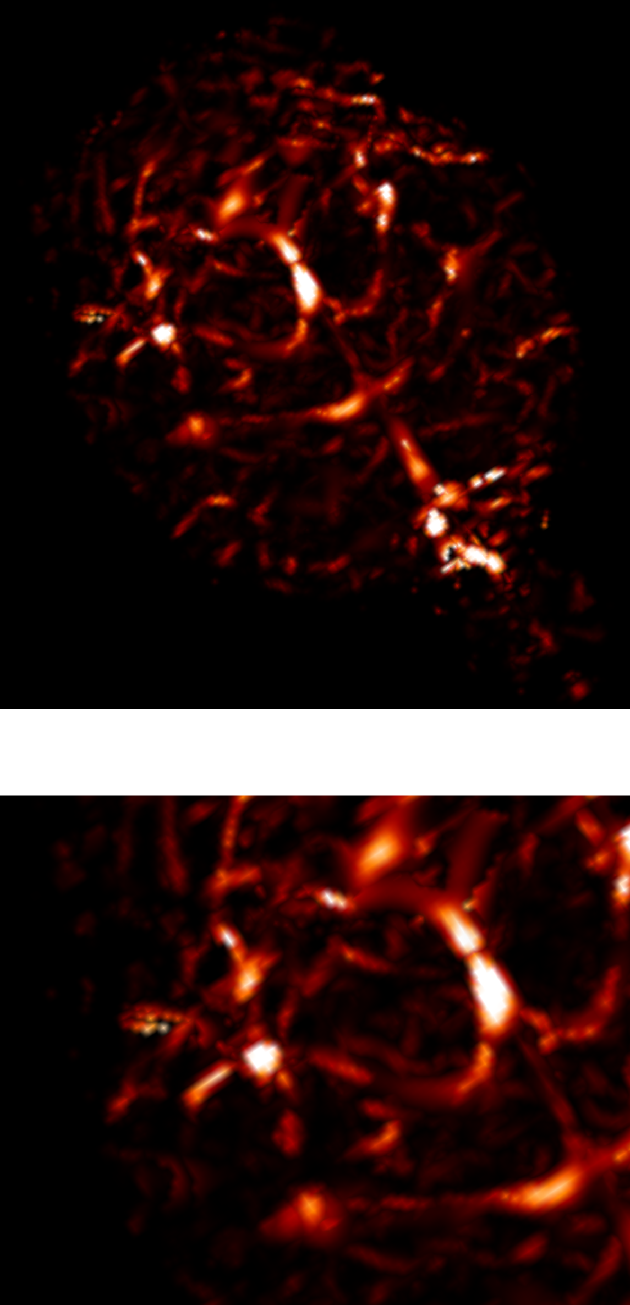}\\
  & PSNR=25.5 dB & PSNR=24.1 dB\\
\end{tabular}\vspace*{-.2cm}
\caption{\label{fig:3dvessels}Angiograms computed from full $k$-space pre-(a) and post-(d) injection data. Angiograms computed from pre-(resp., post-) injection data for decimated $k$-space with $r\!=\!7.3$~(b) and $r\!=\!17.4$~(c)~(resp., (e) and (f)).}
\end{figure}

\vspace*{-.02 \linewidth}
Using the same method as in 2D, namely TSP-sampling and projection onto the set of constraints, we reconstructed volumes from 3D $k$-space. In order to estimate the quality of the reconstructions, we compared the angiograms computed from the 3D images using Frangi filtering~\cite{Frangi98}. The results are shown in Fig.~\ref{fig:3dvessels} for acceleration factors $r=7.3$ (Fig.~\ref{fig:3dvessels}(b,e)) and $r=17.4$ (Fig.~\ref{fig:3dvessels}(c,f)) and compared to the angiogram computed from the whole data.\\
\indent Using the strategy described in Part~\ref{part:hypothesis} the time to traverse $k$-space would be 3.53~s (full acquisition), 3.15~s ($r=7.2$) and 0.88~s ($r=17$). The main drawback of TSP-based sampling schemes is that the time reduction is not directly proportional to $r$, in contrast to classical 2D downsampling and reading out along the third dimension. Nevertheless, if the number of measurements is fixed, the TSP-based approach leads to more accurate reconstruction results since the sampling scheme may fit any density~\cite{Chauffert14}.\\
\indent Angiograms shown in Fig.~\ref{fig:3dvessels} illustrate that one can reduce the travel time in the $k$-space and still observe accurate microvascular structure. If $r=7.3$, time reduction is minor (about 10\% less), but the computed angiogram is almost the same as the one obtained with a complete $k$-space. It is interesting to notice that with a higher acceleration factor~($r=17.4$), the acquisition time is reduced by 75\%, but the computed angiogram remains of good quality. The angiogram appears a bit noisier, especially in the pre-injection setting~(Fig.~\ref{fig:3dvessels}(c)), but the post-injection image allows recovering Willis polygon and most of the major vessels of the mouse brain~(Fig.~\ref{fig:3dvessels}(f)).

\section{CONCLUSION AND FUTURE WORK}
\vspace*{-.03\linewidth}
In this paper, we have shown that the projection algorithm introduced in~\cite{Chauffert14b} is a promising technique to design gradient waveforms. In particular, it allows us to design gradient waveforms in order to implement TSP-based VDS, while existing gradient design methods lead to extremely large acquisition time.
We are currently implementing these waveforms on actual scanners to validate our method on a real CS-MRI framework.
\vspace*{-.02\linewidth}

\section{Acknowledgments}
\vspace*{-.03\linewidth}
This work benefited from the support of the ”FMJH Program Gaspard
Monge in optimization and operation research”, and from EDF's support to
this program.

{\footnotesize
\bibliographystyle{IEEEbib}

}
\end{document}